\documentclass{amsart}

\usepackage{amssymb}

\newtheorem{thm}{Theorem}

\newtheorem{lem}[thm]{Lemma}
\newtheorem{cor}[thm]{Corollary}

\newtheorem{prop}[thm]{Proposition}

\theoremstyle{definition}
\newtheorem{defn}[thm]{Definition}

\newtheorem{say}[thm]{}
\newtheorem{exmp}[thm]{Example}

\newtheorem*{ack}{Acknowledgments}      

\newtheorem{defn-thm}[thm]{Definition--Theorem}  
\newtheorem{defn-lem}[thm]{Definition--Lemma}  

\newtheorem{ass}[thm]{Assumption}

\theoremstyle{remark}

\renewcommand{\c}[0]{{\mathbb C}}  

\renewcommand{\o}[0]{{\mathcal O}} 

\renewcommand{\r}[0]{{\mathbb R}} 

\renewcommand{\a}[0]{{\mathbb A}}

\newcommand{\q}[0]{{\mathbb Q}}
\newcommand{\map}[0]{\dasharrow}
\newcommand{\qtq}[1]{\quad\mbox{#1}\quad}
\newcommand{\spec}[0]{\operatorname{Spec}}

\newcommand{\rank}[0]{\operatorname{rank}}

\newcommand{\supp}[0]{\operatorname{Supp}}    
\newcommand{\red}[0]{\operatorname{red}}    
    
\newcommand{\im}[0]{\operatorname{im}}

\newcommand{\coker}[0]{\operatorname{coker}}

\newcommand{\aut}[0]{\operatorname{Aut}}

\newcommand{\onto}[0]{\twoheadrightarrow}

\newcommand{\rest}[0]{\operatorname{rest}}

\newcommand{\tsum}[0]{\textstyle{\sum}}

\def\into{\DOTSB\lhook\joinrel\to}

\begin{document}
\bibliographystyle{amsalpha}

\hfill\today

\title{Continuous closure of sheaves}
\author{J\'anos Koll\'ar}

\maketitle

\begin{defn}\label{cc.1.defn}
Let $I=(f_1,\dots, f_r)\subset \c[z_1,\dots, z_n]$ be an ideal.
Following  \cite{brenner}
a polynomial $g(z_1,\dots, z_n)$
is in the {\it continuous closure} of $I$ iff there are
 continuous functions $\phi_i$ such that
$g=\phi_1f_1+\cdots+\phi_rf_r$.
These polynomials form an ideal $I^C\supset I$.
For example
$$
z_1^2z_2^2=\frac{\bar z_1 z_2^2}{|z_1|^2+|z_2|^2}z_1^3+
\frac{\bar z_2 z_1^2}{|z_1|^2+|z_2|^2}z_2^3
$$
shows that $z_1^2z_2^2\in (z_1^3,z_2^3)^C\setminus(z_1^3,z_2^3) $.
\end{defn}

The above  definition is very natural, but it is not clear
that it gives an algebraic notion 
(since $\aut(\c/\q)$ does not map continuous functions to
 continuous functions) or that  it defines a sheaf in the Zariski topology
(since a  continuous function may grow faster than any polynomial).

This note has three aims:
\begin{itemize}
\item  We give a purely algebraic construction of
the continuous closure of any torsion free coherent sheaf (\ref{cc.2.defn}).
Although the construction makes sense for any reduced scheme,
even in positive and mixed characteristic,  
it is not clear that it  corresponds to 
a more  intuitive version in  general.

\item  In characteristic 0 we prove that one gets the same definition of $I^C$
 using various subclasses of continuous functions (\ref{mydef=brdef}).

\item We show that taking  continuous  closure  commutes with
flat morphisms whose fibers are  semi-normal (\ref{flat.commute.thm}),
at least in characteristic 0. 
In particular, the  continuous  closure  of a
coherent ideal sheaf is again a coherent ideal sheaf 
 (both in the  Zariski and in the \'etale topologies)
and it commutes with   field extensions.
\end{itemize}

It should be noted that although our definition of the   continuous  closure 
is purely algebraic and without any reference to continuity,
 the proof of these base change properties
uses continuous functions in an essential way.

Instead of working with $\c$ or other algebraically closed fields,
one can also define the   continuous  closure over any
topological field. The most interesting is the real case,
considered in \cite{fef-kol}. The answer turns out to be quite
different; for instance, over $\c$ the  continuous  closure of
$(x^2+y^2)$ is itself but over $\r$ it is the much larger ideal
$(x^2+y^2,x^3, y^3)$. The methods, however, are quite similar.
The main difference is that the base change properties
are not considered in  \cite{fef-kol} and the key construction
(\ref{red.step.main.2}) is  more complicated
over non-closed fields.

The methods of this paper provide a way to compute the
  continuous  closure in principle, but it is unlikely to be
practical in its current form.

\subsection*{Descent problems}{\ }

Instead of working with ideals, I work with
maps of locally free sheaves $f:E\to F$.
Thus an ideal sheaf  $I=(f_1,\dots, f_r)\subset \o_X$
corresponds to the map
$(f_1,\dots, f_r): \o_X^r\to \o_X$.
For inductive purposes we need the case when $E$ and $F$ live on
different schemes.

\begin{defn}\label{desc.prob.efn} Fix a base scheme $S$.
A {\it descent problem}  over $S$ is a compound
object 
$$
{\mathbf D}=\bigl(p:Y\to X, f:p^*E\to F\bigr)
\eqno{(\ref{desc.prob.efn}.1)}
$$ 
 consisting of a 
 proper morphism  $p:Y\to X$ of reduced $S$-schemes of finite type,
 a  locally free sheaf $E$ on $X$,  a
   locally free sheaf $F$ on $Y$
and a   map of  sheaves
$f:p^*E\to F$.

The original setting corresponds to the cases
$$
\bigl(p:X\cong X, f:\o_X^r\to \o_X\bigr)\qtq{with} S=\spec \c
\eqno{(\ref{desc.prob.efn}.2)}
$$
  and, at least when $X$ is (semi)normal,  the continuous closure is 
$$
H^0(X,\o_X)\cap \im\bigl[C^0\bigl(X(\c),\o_X^r\bigr)
\stackrel{f}{\to} C^0\bigl(X(\c),\o_X\bigr)\bigr]
\eqno{(\ref{desc.prob.efn}.3)}
$$
where $C^0$ denotes the space of continuous sections.

Our claim is that in general the primary task should be to understand
the continuous aspects of the problem, that is, 
 the image of
$$
f\circ p^*:C^0\bigl(X(\c),E\bigr)\to C^0\bigl(Y(\c),F\bigr).
\eqno{(\ref{desc.prob.efn}.4)}
$$
Once that is done, the answers to the algebraic questions
should follow.

A  descent problem over $\c$ is called
{\it finitely determined} if for every
$\phi_Y\in C^0\bigl(Y(\c),F\bigr)$ the following are equivalent
\begin{enumerate}
\item[(5.a)] There is a $\phi_X\in C^0\bigl(X(\c),E\bigr)$ such that
$\phi_Y=f\circ p^*(\phi_X)$.
\item[(5.b)] For every finite subset $Z\subset Y$ there is a
$\phi_{X,Z}\in C^0\bigl(X(\c),E\bigr)$ such that
$\phi_Y(z)=f\circ p^*(\phi_{X,Z})(z)$ for every $z\in Z$.
\end{enumerate}

For finitely determined descent problems it is quite easy to pass
between the continuous and the algebraic sides.

The original descent problems (\ref{desc.prob.efn}.2)
are finitely determined only in the trivial case  $I=\o_X$.
A better example is given by the following construction.
Given $I=(f_1,\dots, f_r)$, let $Y:=B_IX$ denote the blow-up of $I$
with projection $p:Y\to X$. The ideal sheaf $f^*I\subset \o_Y$
 is locally free; denote it by $\o_Y(-E)$
where $E$ is an exceptional divisor. We get a descent problem
$$
\bigl(p:Y\to X, f:p^*\o_X^r\to  \o_Y(-E)\bigr),
\eqno{(\ref{desc.prob.efn}.6)}
$$
which is, as we will see,  equivalent to the original one.
 Finite determinacy  for  (\ref{desc.prob.efn}.6) 
is my reformulation of
the axis closure condition of \cite{brenner}
(though they are probably not quite equivalent).
It turns out that  (\ref{desc.prob.efn}.6) is finitely determined
in many cases but not always. Such examples were discovered by \cite{hochster};
an especially nice one is $I=(x^2, y^2, xyz)$.

This paper grew out of first reducing
(\ref{desc.prob.efn}.2) to (\ref{desc.prob.efn}.6)
and then studying the latter by restriction to $E$ and induction.


\end{defn}

The key technical result (\ref{main.thm.v0}) shows that every descent problem
is equivalent to a  finitely determined descent problem.
To achieve this,
we need various ways of modifying descent problems.
The following  definition is chosen to 
consist of simple and computable steps yet
be broad enough
for the proofs to work.
(It should become clear that several variants of the definition would 
also work. The present one is meant to supersede the
choice in \cite{k-obw}.)

\begin{defn}[Scions of  descent problems]\label{rel.desc.prob.sefn}
Let ${\mathbf D}=\bigl(p:Y\to X, f:p^*E\to F\bigr)$  be a descent problem
over $S$.
A {\it scion} of  ${\mathbf D}$
is any  descent problem 
${\mathbf D}_s=\bigl(p_s:Y_s\to X, f_s:p_s^*E\to F_s\bigr)$
that can be obtained by  repeated application of
the following procedures.
\begin{enumerate}
\item   For a proper morphism $r:Y_1\to Y$ set
$$
r^*{\mathbf D}:=
\bigl(p\circ r:Y_1\to X, r^*f: (p\circ r)^*E\to r^*F\bigr).
$$
\item  Given $Y_w$, assume that there are several  proper morphisms 
$r_i:Y_w\to Y$  such that the composites
$p_w:=p\circ r_i$ are all the same. Set
$$
(r_1,\dots, r_m)^*{\mathbf D}:=
\bigl(p_w:Y_w\to X, \tsum_{i=1}^mr_i^*f: p_w^*E\to \tsum_{i=1}^mr_i^*F\bigr)
$$
where $\tsum_{i=1}^mr_i^*f$ is the natural diagonal map.
\item Assume that $f$ factors as 
$p^*E\stackrel{q}{\to} F'\stackrel{j}{\into} F$ 
where 
$F'$ is a locally free sheaf and
$\rank_y j=\rank_y F'$ for all $y$ in  a 
 dense  open subscheme  $Y^0\subset Y$.
Then  set
$$
{\mathbf D}':=\bigl(p:Y\to X, f':=q:p^*E\to F'\bigr).
$$
\end{enumerate}

By construction, each scion
${\mathbf D}_s=\bigl(p_s:Y_s\to X, f_s:p_s^*E\to F_s\bigr)$
comes equipped with a morphism $r_s:Y_s\to Y$, called the
{\it structure map.}

Each scion remembers all of its forebears.
That is, two scions are considered the ``same'' only if they
have been constructed by an identical sequence of procedures.
This is  quite important since the 
way we obtain the
locally free sheaf  $F_s$
 does depend on the
whole sequence.

The class of all scions of ${\mathbf D}$ is denoted by
${\rm Sci}({\mathbf D})$. 

Simple examples of scions are given by {\it restrictions}.
If $Y_1\subset Y$ is a subscheme, we set
$$
{\mathbf D}|_{Y_1}:=\bigl(p|_{Y_1}:Y_1\to X, f|_{Y_1}:
\bigl(p|_{Y_1}\bigr)^*E\to F|_{Y_1}\bigr).
$$
If 
$X_1\subset X$ is a subscheme and $Y_1:=\red p^{-1}(X_1)$, we set
${\mathbf D}|_{X_1}:={\mathbf D}|_{Y_1}$.
\end{defn}

\begin{say}[Seminormalization]\label{semnorm.say}
(For more details, see \cite[Sec.I.7.2]{rc-book}.)
A morphism $p:X'\to X$ 
is a {\it partial seminormalization} if $X'$ is reduced,
$p$ is a finite homeomorphism and 
$k\bigl(p^{-1}(x)\bigr)=k(x)$ for every point $x\in X$.
Under mild conditions (for instance if $X$ is excellent)
there is a unique largest partial seminormalization 
 $\pi:X^{sn}\to X$, called the 
 {\it seminormalization} of $X$.

If $p:Y\to X$ is a proper surjection of reduced schemes
then composing by $p$ identifies $\o_{X^{sn}}$
with those sections of $\o_{Y^{sn}}$ that are constant on the fibers of $p$.

Note that the seminormalization is dominated by the
normalization, thus we can think of the seminormalization
as  a partial normalization.
In some respects, seminormalizations behave better
than the normalization. For instance,
any morphism $g:Y\to X$ induces a morphism
between the seminormalizations $g^{sn}:Y^{sn}\to X^{sn}$.
(For normalization this can fail if $g$ is not dominant.)

A morphism is called seminormal if  its geometric fibers are seminormal.
If $X$ and $g:X'\to X$ are both  seminormal then so is $X'$.
For normal fibers this is proved in \cite[I.7.2.6]{rc-book}.
By localization, the general case follows from 
the following.
\medskip

{\it Claim \ref{semnorm.say}.1.} 
Let $f:(y\in Y)\to (0\in X)$ be a flat morphism of finite type.
Assume that $X$, $Y_0$ and $Y\setminus\{y\}$ are
seminormal. Then $X$ is seminormal.
\medskip

Proof. If $\dim Y_0=0$ then $f$ is smooth and we are done by
\cite[I.7.2.6]{rc-book}. 
Let $h\in \o_{Y^{sn}}$ be a section. 
If $\dim Y_0\geq 1$, we prove by induction
on $r$ that  $h\in \o_Y+m_{0,X}^r\o_{Y^{sn}}$ for every $r$.
We can start with $r=0$. In general, assume that we have
$p_r\in \o_Y$ such that
$h-p_r\in m_{0,X}^r\o_{Y^{sn}}$.
By restricting to $Y_0$ we see that
$(h-p_r)|_{Y_0}$ is a section of
$$
\o_{Y_0^{sn}}\otimes \bigl(m_{0,X}^r/m_{0,X}^{r+1}\bigr)=
\o_{Y_0}\otimes \bigl( m_{0,X}^r/m_{0,X}^{r+1}\bigr).
$$
Thus there is a $q_{r+1}\in m_{0,X}^r\o_Y$
such that $h-p_r-q_{r+1}$ vanishes along $Y_0$ to order $r+1$.
This shows that the completion of $\o_Y$ equals 
  the completion of $\o_{Y^{sn}}$, hence $\o_Y=\o_{Y^{sn}}$.\qed
\medskip

If $F$ is a coherent sheaf on $X$, its pull-back to $X^{sn}$
is denoted by $F^{sn}$. We frequently view $F^{sn}$
as an $\o_X$-sheaf.

If $X$ is a variety over $\c$, then $\o_{X^{sn}}$
consists of those rational functions that are
continuous. Thus it appears that 
the continuous closure is a concept that naturally lives
on seminormal schemes.

It would be  possible to consider descent problems only
for seminormal schemes. This, however, would be inconvenient
since various constructions do not yield
seminormal schemes, and we would have to take
seminormalizations all the time.
Instead, next  we build the seminormalizations into the
definition of the global sections of
${\rm Sci}({\mathbf D})$.
\end{say}

\begin{defn} \label{scion.sect.defn}
Let ${\mathbf D}=\bigl(p:Y\to X, f:p^*E\to F\bigr)$  be a descent problem
with scions
$$
{\rm Sci}({\mathbf D})=
\Bigl\{\bigl(p_i:Y_i\to X, f_i:p_i^*E\to F_i\bigr): i\in I\Bigr\}.
$$
An {\it algebraic   global section} of $F$ over ${\rm Sci}({\mathbf D})$
is a collection of sections
$$
\Phi:=\bigl\{\phi_i\in H^0(Y_i^{sn}, F_i^{sn}): i\in I\bigr\}
$$
 such that
the $\phi_i$ commute with pull-backs for the operations
(\ref{rel.desc.prob.sefn}.1--2)
and with push-forward for the operations
(\ref{rel.desc.prob.sefn}.3).
All sections form an $\o_S$-module
 $$
H^0\bigl({\rm Sci}({\mathbf D}), F\bigr);
$$
one can also think of it as the direct limit
of the $ H^0(Y_i^{sn}, F_i^{sn})$ over the category
${\rm Sci}({\mathbf D})$.
We call $\phi_i$ the {\it restriction} of $\Phi$ to $Y_i$,
denoted by $\Phi|_{Y_i}$.
The most important of these restrictions is $\Phi|_Y$.
Note that $\Phi|_Y$ uniquely determines $\Phi$. 
Indeed, the constructions (\ref{rel.desc.prob.sefn}.1--2) automatically carry
along $\phi$ and in (\ref{rel.desc.prob.sefn}.3) 
the natural map
 $H^0(Y^{sn},F')\to H^0(Y^{sn},F)$ is an injection.

We usually think of $H^0\bigl({\rm Sci}({\mathbf D}), F\bigr)$
as an $\o_S$-submodule of $H^0(Y^{sn}, F^{sn})$.

Note also that every $\phi_X\in H^0(X, E)$
defines a  global section of $F$ over ${\rm Sci}({\mathbf D})$
by setting $\phi_i:=f_i(p_i^*\phi_X)$. Thus we have  natural maps
$$
H^0(X, E)\to H^0\bigl({\rm Sci}({\mathbf D}), F\bigr)
\into H^0(Y^{sn}, F^{sn}).
\eqno{(\ref{scion.sect.defn}.1)}
$$

\end{defn}

We can now define a notion of
continuous closure of sheaves. A justification of the
definition will be given only later in (\ref{mydef=brdef}).

\begin{defn}[Continuous closure of sheaves]\label{cc.2.defn}
Let $X$ be a pure dimensional, reduced, affine scheme 
 over a field of characteristic 0 and $J$ a torsion free 
coherent sheaf on $X$. One can realize $J$ as the
image of a map between locally free shaves $f:E\to F$.
Let ${\mathbf D}_J=\bigl(p:Y \cong X, f:E\to F\bigr)$ be the
corresponding  descent problem.
Define the {\it continuous closure} of $J$ as
$$
J^C:=H^0\bigl({\rm Sci}({\mathbf D}_J), F\bigr)
\subset H^0\bigl(X^{sn}, F^{sn}\bigr).
$$
We see later (\ref{cl.indep.say}) that $J^C$ does not depend on the choice
of $f:E\to F$.

\end{defn}

The above definition is purely algebraic but it does not
connect with continuity in any obvious way.
Actually, for base fields that are not naturally
subfields of $\c$, it is not even clear what continuity should mean.
This is the question we consider next.

\subsection*{Classes of  continuous functions}{\ }

Here we describe various classes of  functions
where out proof works.

\begin{ass}\label{function.classes.ass}
Let $k$ be a field and $K\supset k$ an algebraically closed field.
For a $k$-scheme of finite type, let $C^K(X)$ denote the 
$K$-vector space 
of all functions $X(K)\to K$.
We consider vector subspaces
$C^*\bigl(Z\bigr)\subset C^K\bigl(Z\bigr)$
 that satisfy
the following properties.
\begin{enumerate}
\item (Sheaf) If $Z=\cup_i U_i$ is an open cover of $Z$
then $\phi\in C^*\bigl(Z\bigr)$ iff
$\phi|_{U_i}\in C^*\bigl(U_i\bigr)$ for every $i$.
\item ($\o_Z$-module) If $\phi\in C^*\bigl(Z\bigr)$
and $h\in \o_Z$ is a regular function  then
 $h\cdot \phi\in C^*\bigl(Z\bigr)$.
\item (Pull-back) For every $k$-morphism 
$g:Z_1\to Z_2$, composing with 
$g$ maps
$C^*\bigl(Z_2\bigr)$ to $C^*\bigl(Z_1\bigr)$.
\item (Zariski dense is dense)
  Let $\phi\in C^*\bigl(Z\bigr)$ and $h$ a rational function
on $Z$ such that $\phi$ equals $h$ on a  dense open subset. Then
$\phi=h$ everywhere and $h$ is a regular function on $Z^{sn}$.
This also implies that the support of  every $\phi\in C^*\bigl(Z\bigr)$
 is a union of irreducible components of $Z$.
\item (Descent property)
Let $g:Z_1\to Z_2$ be a proper,  dominant   $k$-morphism, 
 $\phi\in C^K\bigl(Z_2\bigr)$ and assume that
$\phi\circ g\in C^*\bigl(Z_1\bigr)$.
Then $\phi\in C^*\bigl(Z_2\bigr)$.

In particular, assume that $X$ is a union of its closed subvarieties
$X_i$ and we have $\phi_i\in C^*(X_i)$ such that 
$\phi_i|_{X_i\cap X_j}=\phi_j|_{X_i\cap X_j}$ for every $i,j$.
The descent property for $\amalg_i X_i\to X$ shows that
there is a $\phi\in C^*(X)$ such that 
$\phi|_{X_i}=\phi_i$ for every $i$.
\item (Extension property) Let $Z_1\subset Z_2$ be a closed 
 subscheme.
Then the restriction map
$C^*\bigl(Z_2\bigr)\to C^*\bigl(Z_1\bigr) $ is surjective.
\item (Cartan--Serre A and B) Every locally free sheaf is generated by
finitely many $C^*$-sections and every surjection of
locally free sheaves has a $C^*$-valued splitting.
(For more details, see (\ref{basic.props}).)
\end{enumerate}
We can unite (5) and (6)  as follows.
\begin{enumerate}
\item[(5+6)] (Strong descent property) 
Let $g:Z_1\to Z_2$ be a proper   $k$-morphism and 
 $\psi\in C^*\bigl(Z_2\bigr)$. Then
$\psi=\phi\circ g$ for some  $\phi\in C^*\bigl(Z_2\bigr)$
iff $\psi$ is constant on every fiber of $g$.
\end{enumerate}
\end{ass}

\begin{exmp}\label{function.classes.exmp} 
Here are some natural examples satisfying the
assumptions (\ref{function.classes.ass}.1--7).
Let us start with the cases when $k\subset K=\c$.
\begin{enumerate}
\item  Let
$C^0\bigl(Z\bigr)$ denote all
continuous functions on  $Z(\c)$.
\item  Let $C^h\bigl(Z\bigr)$ denote all locally 
H\"older continuous functions on $Z(\c)$.
\item Let
$S^0(X)$ be the sheaf of $\c$-valued  continuous 
semi-algebraic functions  on
 $X(\c)$, viewed as a real algebraic variety.  
(If $X\subset \c^m$, we identify $\c^m$ with
$\r^{2m}$ and view $X(\c)$ as a real variety.
A function on $\r^{2m}$ is semi-algebraic iff its graph
is semi-algebraic, that is, a finite union of
sets defined by polynomial inequalities of the form $f\geq 0$.) 
See  \cite[Chap.2]{bcr} for details and proofs of the properties
(\ref{function.classes.ass}.1--7).
(Let me just note that (\ref{function.classes.ass}.4)
is more interesting than it sounds.
For instance, on the Whitney umbrella  $(x^2=y^2z)\subset \r^3$
not every Zariski dense open set is Euclidean dense.)
\end{enumerate}
I do not know how to generalize the first two of these
in case $k$ is not embedded into $\c$, but the third variant
can be extended to any characteristic 0 field.
\begin{enumerate}\setcounter{enumi}{3}
\item Let ${\mathbf R}$ be a real closed field,
 ${\mathbf C}:={\mathbf R}\bigl(\sqrt{-1}\bigr)$
and assume that $k\subset  {\mathbf C}$.  Let
$S^0_{\mathbf R}(X)$ be the sheaf of ${\mathbf C}$-valued  continuous 
semi-algebraic functions  on
 $X({\mathbf C})$, viewed as an
 ${\mathbf R}$-variety.  (See  \cite[Chap.2]{bcr} for details.)
\end{enumerate}
I do not know any examples in positive characteristic.
\end{exmp}


\begin{say}[$C^*$-valued sections]\label{basic.props}
Let $F$ be a  locally free sheaf on $Z$
and $Z=\cup_i U_i$ an open cover such that
$F|_{U_i}$ is trivial of rank $r$ for every $i$. Let
$C^*\bigl(Z, F\bigr)$
denote the set of those sections 
such that
$\phi|_{U_i}\in C^*\bigl(U_i\bigr)^r$ for every $i$.
If $C^*$ satisfies the properties (\ref{function.classes.ass}.1--2),
 this is independent of the
trivializations and the choice of the covering.

Assume next that 
(\ref{function.classes.ass}.7) holds.
We claim that if
  $C^*$ satisfies the properties (\ref{function.classes.ass}.1--6)
then their natural analogs also hold for $C^*\bigl(Z, F\bigr)$.
This is clear for the properties (\ref{function.classes.ass}.2--5).

In order to check the extension property (\ref{function.classes.ass}.6),
  let  $Z_1\subset Z_2$ be an closed subvariety
and $F$ a locally free sheaf on $Z_2$. Write it as  a quotient of a
trivial bundle $\o_{Z_2}^N$. Every section 
$\phi_1\in C^*\bigl(Z_1, F|_{Z_1}\bigr)$ lifts to a
section in $C^*\bigl(Z_1, \o_{Z_1}^N\bigr)$
which in turn extends to  a
section in $C^*\bigl(Z_2, \o_{Z_2}^N\bigr)$ by (\ref{function.classes.ass}.6).
The image of this lift in $C^*\bigl(Z_2, F|_{Z_2}\bigr)$
gives the required lifting of  $\phi_1$.

Let ${\mathbf D}$ be a descent problem with
scions ${\rm Sci}({\mathbf D})$.
If $\phi\in C^*(Y,F)$ then $r^*\phi\in C^*(Y_1,r^*F)$
and $\tsum_{i=1}^mr_i^*\phi\in C^*(Y_w,\tsum_{i=1}^mr_i^*F)$
are well defined. 
In (\ref{rel.desc.prob.sefn}.3) above,  
$j:C^*(Y,F')\to C^*(Y,F)$ is an injection,
hence there is at most one $\phi'\in C^*(Y,F')$
such that $j(\phi')=\phi$.
Iterating these, for any scion ${\mathbf D}_s$ of ${\mathbf D}$
we get a partially defined map, called the {\it restriction},
$$
\rest: C^*(Y,F)\map C^*(Y_s, F_s)\qtq{denoted by}
\phi\mapsto \phi|_{Y_s}\qtq{or} \phi\mapsto \phi|_{{\mathbf D}_s}.
$$

The restriction map sits in a commutative square
$$
\begin{array}{ccc}
C^*(Y,F) & \stackrel{\rest}{\map} & C^*(Y_s, F_s)\\
\uparrow && \uparrow \\
C^*(X,E) & = & C^*(X,E).
\end{array}
$$
If the structure map $r_s:Y_s\to Y$ is surjective then
 the restriction map
$\rest: C^*(Y,F)\map C^*(Y_s, F_s)$ is injective (on its domain).
In this case,  understanding the image of
$f\circ p^*:C^*(X,E)\to C^*(Y,F)$ is pretty much
equivalent to understanding the image of
$f_s\circ p_s^*:C^*(X,E)\to C^*(Y_s,F_s)$.

As  long as $C^*$ satisfies  the properties 
(\ref{function.classes.ass}.1--3),
we can follow the definition
(\ref{scion.sect.defn}) to  obtain
$$
C^*\bigl({\rm Sci}({\mathbf D}), F\bigr),
\eqno{(\ref{basic.props}.1)}
$$
the space of $C^*$-valued global sections of $F$ over ${\rm Sci}({\mathbf D})$.
 We have  natural maps
$$
C^*(X, E)\to C^*\bigl({\rm Sci}({\mathbf D}), F\bigr)
\into C^*(Y^{sn}, F^{sn})\ \bigl(\ = C^*(Y, F)\bigr).
\eqno{(\ref{basic.props}.2)}
$$
Note further that
$$
H^0\bigl({\rm Sci}({\mathbf D}), F\bigr)=
C^*\bigl({\rm Sci}({\mathbf D}), F\bigr)
\cap H^0(Y^{sn}, F^{sn}).
\eqno{(\ref{basic.props}.3)}
$$
To see this we need to show that if
$\Phi\in C^*\bigl({\rm Sci}({\mathbf D}), F\bigr)$
and $\Phi|_Y$ is algebraic then every other restriction of $\Phi$
is also algebraic. This is clear for the steps
(\ref{rel.desc.prob.sefn}.1--2). 
For scions as in (\ref{rel.desc.prob.sefn}.3), let
$\phi$ be an algebraic section of $F$. We assume that
$\phi$ is a $C^*$-valued section of $F'$. It is also
a rational section over a Zariski dense open set, thus, by
(\ref{function.classes.ass}.4)
$\phi$ is also an algebraic section of $F'$.

The restriction map on $C^*(Y,F)$ gives a restriction map
on global sections of scions
which also sits in a commutative diagram
$$
\begin{array}{ccc}
C^*\bigl({\rm Sci}({\mathbf D}), F\bigr) & \stackrel{\rest}{\to} & 
C^*\bigl({\rm Sci}({\mathbf D}_s), F_s\bigr)\\
\uparrow && \uparrow \\
C^*(X,E) & = & C^*(X,E).
\end{array}
\eqno{(\ref{basic.props}.4)}
$$ 
Note that the restriction map
on global sections of scions is everywhere defined.
In essence, we defined $ C^*\bigl({\rm Sci}({\mathbf D}), F\bigr)$
to ensure this.

\end{say}

\subsection*{Finitely determined descent problems}{\ }

The notion of a
finitely determined descent problem
(\ref{desc.prob.efn}.5) admits an obvious generalization
to the $C^*$-valued case.
We also need the following more general version.

\begin{defn}\label{findet.defn}  
Let ${\mathbf D}=\bigl(p:Y\to X, f:p^*E\to F\bigr)$
 be a  descent problem and $Z\subset X$ a 
 closed algebraic subvariety. ${\mathbf D}$ is called
{\it finitely determined} relative to $Z$  if for every
$\phi_Y\in C^*\bigl(Y,F\bigr)$ 
that vanishes on $p^{-1}(Z)$ the following are equivalent
\begin{enumerate}
\item There is a $\phi_X\in C^*\bigl(X,E\bigr)$ 
such that $\phi_Y=f\circ p^*(\phi_X)$.
\item For every finite subset $\{y_1,\dots, y_m\}\subset Y$ there is a
$\phi_{X,y_1,\dots, y_m}\in C^*\bigl(X,E\bigr)$ such that
$\phi_Y(y_i)=f\circ p^*(\phi_{X,y_1,\dots, y_m})(y_i)$ for  $i=1,\dots, m$.
\end{enumerate}
We see in (\ref{wronsk.lem}) that these are also
equivalent to the following precise form:
\begin{enumerate}\setcounter{enumi}{2}
\item The above (2) holds for all $m\leq \rank E+1$.
\end{enumerate}
\end{defn}

Although  (\ref{findet.defn}.2)
asks about all possible finite sets of points in $Y$,
the conditions imposed by points in different fibers
of $p$ are independent. Thus the only interesting case
is when all the $y_i$ are in the same fiber.
Working in a fiber, we have a general abstract test.

\begin{lem}[Wronskian test]\label{wronsk.lem} Let $Y$ be a set and
 $\phi,f_1,\dots, f_r$   functions on  $Y$ with values in a field $K$.
Assume that the $f_i$ are linearly independent.
Then the following are equivalent.
\begin{enumerate}
\item  $\phi$ is a  linear combination of the $f_i$.
\item For every $r+1$ points  $y_1,\dots, y_{r+1}$ there are
$c_1,\dots, c_r$ (possibly depending on the $y_i$)
such that
$\phi(y_i)=\sum_j c_jf_j(y_i)$ for $i=1,\dots, r+1$.
\item The following  determinant is identically zero as a function on $Y^{r+1}$.
$$
\left|
\begin{array}{cccc}
f_1(y_1) & \cdots & f_1(y_r) &  f_1(y_{r+1})\\
\vdots &&\vdots & \vdots \\
f_r(y_1) & \cdots & f_r(y_r) &   f_r(y_{r+1})\\
\phi(y_1) &\cdots & \phi(y_r) &   \phi(y_{r+1})
\end{array}
\right|
$$
\end{enumerate}
\end{lem}

Proof. Since the $f_i$ are linearly independent, 
there are  $y_1,\dots, y_r\in Y$
such that 
the upper left $r\times r$ subdeterminant above is nonzero.
Fix these  $y_1,\dots, y_r$ and solve the linear system
$$
\phi(y_i)=\tsum_j\ 
\lambda_j f_j(y_i)\qtq{for $i=1,\dots,r$.}
$$
Replace $\phi$ by $\psi:=\phi-\sum_i \lambda_if_i$
and let $y_{r+1}$ vary.
Then our determinant is
$$
\left|
\begin{array}{cccc}
f_1(y_1) & \cdots & f_1(y_r) &  f_1(y_{r+1})\\
\vdots &&\vdots & \vdots \\
f_r(y_1) & \cdots & f_r(y_r) &   f_r(y_{r+1})\\
0 &\cdots & 0 &   \psi(y_{r+1})
\end{array}
\right|
$$
The whole determinant vanishes iff $\psi(y_{r+1})$
is identically zero.
That is, when $\phi\equiv \sum_j \lambda_j f_j$.\qed
\medskip

If a  descent problem is not finitely determined,
we can still study the conditions imposed by
(\ref{findet.defn}.2). This leads to the
following definition.

\begin{defn}\label{sci.0.defn}
Given a  descent problem  ${\mathbf D}=\bigl(p:Y\to X, f:p^*E\to F\bigr)$,
let ${\rm Sci}^0({\mathbf D})\subset {\rm Sci}({\mathbf D})$
denote all $0$-dimensional scions and  ${\mathbf D}$ itself.
We can now define
$$
H^0\bigl({\rm Sci}^0({\mathbf D}), F\bigr)\qtq{and}
C^*\bigl({\rm Sci}^0({\mathbf D}), F\bigr)
$$
as the collection of sections
$\bigl\{\phi_i\in H^0\bigl(Y_i^{sn}, F_i^{sn}\bigr)\bigr\}$
(resp.\ $\bigl\{\phi_i\in C^*\bigl(Y_i, F_i\bigr)\bigr\}$)
that satisfy the compatibility conditions
as in (\ref{scion.sect.defn})
where now $Y_i$ runs through only the scions in
${\rm Sci}^0({\mathbf D})$.

Thus  ${\mathbf D}$ is finitely determined iff
$$
\im\bigl[C^*\bigl(X,E\bigr)\to C^*\bigl(Y,F\bigr)\bigr]=
C^*\bigl({\rm Sci}^0({\mathbf D}), F\bigr).
$$
\end{defn}

An advantage of $H^0\bigl({\rm Sci}^0({\mathbf D}), F\bigr)$
is that it can be easily computed algebraically.

\begin{say}[Computation of $H^0\bigl({\rm Sci}^0({\mathbf D}), F\bigr)$]
\label{exact.seq.say}
Let $X$ be an affine scheme
of finite type  over a field  and 
  ${\mathbf D}=\bigl(p:Y\to X, f:p^*E\to F\bigr)$  a 
descent problem.
 We inductively construct descent problems
${\mathbf D}_i=\bigl(p_i:Y_i \to X_i, f_i:p_i^*E_i\to F_i\bigr)$
as follows. Set
${\mathbf D}_0:={\mathbf D}$ and assume that 
${\mathbf D}_i$ is already constructed.

By Cohomology and Base Change  \cite[III.12.11]{hartsh}
there is a largest open dense subset $X_i^0\subset X_i$
 over  which the following hold:
\begin{enumerate}
\item $p_i^{sn}:Y_i^{sn}\to X_i$ is flat,
\item the $R^j(p_i^{sn})_*F_i^{sn}$ are locally free and 
commute with base change, and 
\item $E_i\to (p_i^{sn})_*F_i^{sn}$ has constant rank.
\end{enumerate}
Set $X_{i+1}:=X_i\setminus X_i^0$ and let
${\mathbf D}_{i+1}$ be the restriction of ${\mathbf D}_i$
to $X_{i+1}$.

Set $Q_i^0:=\coker\bigl[E_i\to (p_i^{sn})_*F_i^{sn}\bigr]$
and let $Q_i$ be the push forward of  $Q_i^0$
by the locally closed embedding
$ X_i^0\into X$. The $Q_i$ are quasi-coherent sheaves on $X$.
We get natural sheaf maps
$q_i:E\to E_i\to Q_i$.

By construction, 
if $x\in X_i^0$ and $\phi\in H^0(Y^{sn},F^{sn})$
then $\phi$ satisfies (\ref{findet.defn}.2) for all
subsets of $p^{-1}(x)$ iff $q_i(\phi)\in H^0(X, Q_i)$
vanishes at $x$. Since $X=\cup_iX_i^0$,
this implies that
$$
H^0\bigl({\rm Sci}^0({\mathbf D}), F\bigr)= \ker
\bigl[ H^0\bigl(Y^{sn},F^{sn}\bigr)\to \tsum_i H^0\bigl(X,Q_i\bigr)\bigr].
\eqno{(\ref{exact.seq.say}.4)}
$$
\end{say}

This implies important functoriality properties of
$H^0\bigl({\rm Sci}^0({\mathbf D}), F\bigr)$,
but first we need a definition.

\begin{defn}[Pulling back descent problems] 
Let ${\mathbf D}=\bigl(p:Y\to X, f:p^*E\to F\bigr)$ be a 
descent problem over a base field $k$.
We consider two ways of obtaining new  descent problems
by base change.

 First, every field extension
$k'\supset k$ gives a descent problem over  $k'$
$$
{\mathbf D}_{k'}:=\bigl(p_{k'}:Y_{k'}\to X_{k'}, 
f_{k'}:p_{k'}^*E_{k'}\to F_{k'}\bigr).
$$

Second, let $b:X^0\to X$ be a
  flat, finite type  morphism  with reduced fibers. 
Let $Y$ be a reduced scheme and $p:Y\to X$ a morphism.
Then $b_Y: Y':=X'\times_XY\to Y$ is flat with reduced fibers,
hence $Y'$ is also reduced.
Thus
$$
b^*{\mathbf D}:=\bigl(p':Y'\to X', f':(p')^*E\to b_Y^*F\bigr)
$$
is also a descent problem.
All the constructions in (\ref{rel.desc.prob.sefn}) commute
with pull-back by flat morphisms  with reduced fibers.
Thus we get a pull-back map
$b^*:{\rm Sci}({\mathbf D})\to {\rm Sci}\bigl(b^*{\mathbf D}\bigr)$.

Note that it is not obvious that there is a
pull-back map
$b^*:H^0\bigl({\rm Sci}({\mathbf D}), F\bigr)\to 
H^0\bigl({\rm Sci}(b^*{\mathbf D}), b_Y^*F\bigr)$.
(Indeed,  $b^*{\mathbf D}$ has scions that are not pulled-back from
${\rm Sci}({\mathbf D})$ and these could pose additional
restrictions on sections.)
We see in (\ref{faithful=isom}) that such problems do not arise.
\end{defn}

\begin{prop}\label{flat.commute.cor}
Let $X$ be an affine scheme
of finite type  over a field and 
  ${\mathbf D}=\bigl(p:Y\to X, f:p^*E\to F\bigr)$  a 
descent problem.
Then the formation of $H^0\bigl({\rm Sci}^0({\mathbf D}), F\bigr)$
commutes with flat, seminormal base changes
and with base field extensions.
\end{prop}

Proof. Note that in (\ref{exact.seq.say})
the formation of the $X_i$ and $Q_i$
 commutes with flat, seminormal base changes
and with base field extensions.
Using (\ref{exact.seq.say}.4), this implies that
 $H^0\bigl({\rm Sci}^0({\mathbf D}), F\bigr)$
also  commutes with flat, seminormal base changes
and with base field extensions. \qed
\medskip

\subsection*{The main theorem and its consequences}{\ }

\begin{defn}[Universal properties]
Let ${\mathbf P}$ be a property of descent problems.
Let ${\mathbf D}$ be a  descent problem over a field $k$.
We say that ${\mathbf D}$ is {\it universally} ${\mathbf P}$
if $b^*{\mathbf D}_{k'}$ satisfies  ${\mathbf P}$ for every 
 base field extension $k'\supset k$
followed by any
 flat, finite type, seminormal base change $b:X'_{k'}\to X_{k'}$.
\end{defn}

The  main technical result of this note is the following.

\begin{thm}\label{main.thm.v0}
 Let ${\mathbf D}=\bigl(p:Y\to X, f:p^*E\to F\bigr)$
be a descent problem over a field of characteristic 0.  
Then it has a universally finitely determined scion
${\mathbf D}_s=\bigl(p_s:Y_s\to X, f_s:p_s^*E\to F_s\bigr)$
whose structure map $r_s:Y_s\to Y$ is surjective.
\end{thm}

Before giving a proof, let us consider some
consequences. 
First we have 
the following property, which was the 
very reason for our definition of scions.

\begin{cor}\label{main.thm.v1}
 Let ${\mathbf D}=\bigl(p:Y\to X, f:p^*E\to F\bigr)$
be a descent problem over a field of characteristic 0.  Assume that 
$C^*$ satisfies  the properties (\ref{function.classes.ass}.1--7).
Then
$$
 C^*\bigl({\rm Sci}({\mathbf D}), F\bigr)=
\im\bigl[C^*(X, E)\to C^*(Y, F)\bigr].
$$
\end{cor}

Proof. Note that, by (\ref{basic.props}.2), the containment
$$
 C^*\bigl({\rm Sci}({\mathbf D}), F\bigr)\supset
\im\bigl[C^*(X, E)\to C^*(Y, F)\bigr]
$$
always holds. 
To see the converse, 
let ${\mathbf D}_s=\bigl(p_s:Y_s\to X, f_s:p_s^*E\to F_s\bigr)$
be a   finitely determined scion of
${\mathbf D}$ whose structure map $r_s:Y_s\to Y$ is surjective. 
We have the obvious inclusions
$$
C^*\bigl({\rm Sci}({\mathbf D}), F\bigr)\subset
C^*\bigl({\rm Sci}({\mathbf D}_s), F_s\bigr)
\subset C^*\bigl({\rm Sci}^0({\mathbf D}_s), F_s\bigr)
$$
and 
$$
 C^*\bigl({\rm Sci}^0({\mathbf D}_s), F_s\bigr)=
\im\bigl[C^*(X, E)\to C^*(Y_s, F_s)\bigr]
$$
since ${\mathbf D}_s$ is finitely determined. 
Note further that $C^*(X, E)\to C^*(Y_s, F_s)$
factors through $C^*(Y, F)$ and through
$C^*\bigl({\rm Sci}({\mathbf D}), F\bigr)$.
Since  the structure map $r_s:Y_s\to Y$ is surjective, 
$C^*(Y, F)\to C^*(Y_s, F_s)$ is injective. These show that
$$
 C^*\bigl({\rm Sci}({\mathbf D}), F\bigr)\subset
\im\bigl[C^*(X, E)\to C^*(Y, F)\bigr]. \qed
$$

We can now see that the 2 definitions of the continuous closure,
(\ref{cc.2.defn}) and the obvious generalization of 
(\ref{cc.1.defn}),   agree with each other.

\begin{cor}\label{mydef=brdef}  Let $X$ be a reduced affine scheme 
 over a field of characteristic 0 and
$f:E\to F$ a map between locally free sheaves.
Set $J=\im(f)$, as a subsheaf of $F$. Then
$$
J^C=\im\bigl[C^*(X, E)\to
 C^*(X, F)\bigr]\cap H^0(X^{sn}, F^{sn}).
$$
\end{cor}

Proof. By definition,
 $J^C=H^0\bigl({\rm Sci}({\mathbf D}), F\bigr)$ and, by (\ref{basic.props}.3), 
$$
H^0\bigl({\rm Sci}({\mathbf D}), F\bigr)=
C^*\bigl({\rm Sci}({\mathbf D}), F\bigr)\cap H^0(X^{sn}, F^{sn}).
$$
By  (\ref{main.thm.v1}),
$C^*\bigl({\rm Sci}({\mathbf D}), F\bigr)=
\im\bigl[C^*(X, E)\to C^*(X, F)\bigr]$. \qed

\medskip

As another consequence, we obtain  that
 global sections of scions are unchanged by  surjective  structure maps. 
Note that we use the invariance of continuous sections to
derive the invariance of algebraic sections.

\begin{cor}\label{faithful=isom}  

Let ${\mathbf D}$ be a descent problem over an affine
base $X$ over a field of characteristic 0.
Let ${\mathbf D}_s$ be a  scion of
${\mathbf D}$ whose structure map $r_s:Y_s\to Y$ is surjective.
 Then the restriction maps
$$
 C^*\bigl({\rm Sci}({\mathbf D}), F\bigr)\to
C^*\bigl({\rm Sci}({\mathbf D}_s), F_s\bigr)
\qtq{and}
 H^0\bigl({\rm Sci}({\mathbf D}), F\bigr)\to
H^0\bigl({\rm Sci}({\mathbf D}_s), F_s\bigr)
$$
are isomorphism.
\end{cor}

Proof. Since $r_s$ is surjective, the restriction maps are
 injective.  By (\ref{main.thm.v1}),
$C^*(X,E)\to C^*\bigl({\rm Sci}({\mathbf D}_s), F_s\bigr)$ is surjective
and it factors through $C^*\bigl({\rm Sci}({\mathbf D}), F\bigr)$.
Thus the restriction map is surjective with $C^*$-coefficients.

The algebraic case also follows once we prove that if
$\phi\in  C^*\bigl({\rm Sci}({\mathbf D}), F\bigr)$
and its restriction to $Y_s$ is algebraic then
$\phi$ itself is algebraic.  This is  a local question on $Y$,
hence we need to show that if $\phi\in C^*(Y)$ and
$r_s^*\phi$ is  a regular function then $\phi$ is a regular function
on $Y^{sn}$.
We can view $\phi$ as a morphism to $\a^1_Y$; let $Y'$ be its image.
Since $Y_s\to Y$ is proper, $Y'\to Y$ is  proper 
and $Y'_K\to Y_K$ is a homeomorphism.
 Thus $Y'$ is dominated by the
seminormalization.
\qed
\medskip

The next result is an important invariance property
of global sections of descent problems.

\begin{cor} \label{flat.commute.thm}
Let $X$ be an affine scheme
of finite type  over a field of characteristic 0 and 
  ${\mathbf D}=\bigl(p:Y\to X, f:p^*E\to F\bigr)$  a 
descent problem.
Then taking algebraic global sections of ${\rm Sci}({\mathbf D})$
commutes  with base field extensions and with flat, seminormal base changes.

In particular, taking the continuous closure 
commutes  with base field extensions and with flat, seminormal base changes.
\end{cor}

Proof. By (\ref{main.thm.v0}),
  ${\mathbf D}$ has a universally finitely determined scion
${\mathbf D}_s=\bigl(p_s:Y_s\to X, f_s:p_s^*E\to F_s\bigr)$
whose structure map $r_s:Y_s\to Y$ is surjective.

By (\ref{faithful=isom}),
$$
H^0\bigl({\rm Sci}({\mathbf D}), F\bigr)=
H^0\bigl({\rm Sci}({\mathbf D}_s), F_s\bigr)
$$
and the equality continues to hold after every base change.
Thus it is sufficient to prove (\ref{flat.commute.thm})
in case  ${\mathbf D}$ is universally finitely determined.
For such descent problems
$$
H^0\bigl({\rm Sci}({\mathbf D}), F\bigr)=
H^0\bigl({\rm Sci}^0({\mathbf D}), F\bigr),
$$
and we saw in (\ref{flat.commute.cor})
 that
$H^0\bigl({\rm Sci}^0({\mathbf D}), F\bigr)$
 commutes  with base field extensions and with flat, seminormal base changes.
\qed
\medskip

Since open embeddings are flat with seminormal fibers,
we can sheafify the notion of continuous closure.

\begin{defn}\label{flat.commute.def}
Let ${\mathbf D}=\bigl(p:Y \to X, f:p^*E\to F\bigr)$ be a
  descent problem.
By  (\ref{flat.commute.thm}),
 as $\{U:U\into X\}$ runs through all
 affine open subsets,
the rule
$$
U\mapsto H^0\bigl({\rm Sci}({\mathbf D}|_U), F|_U\bigr)
$$
defines a coherent sheaf in the Zariski topology,  denoted by
$$
R^0p_*\bigl({\rm Sci}({\mathbf D}), F\bigr)
\eqno{(\ref{flat.commute.def}.1)}
$$
and called the {\it push forward} of ${\rm Sci}({\mathbf D})$.

As in (\ref{scion.sect.defn}.1), there are natural maps
$$
E\to R^0p_*\bigl({\rm Sci}({\mathbf D}), F\bigr)
\into p_*F^{sn}.
\eqno{(\ref{flat.commute.def}.2)}
$$
\end{defn}

Finally, let us see that the definition
(\ref{cc.2.defn}) is independent of the auxiliary choices.

\begin{prop}\label{cl.indep.say}
The continuous closure is independent of the choice of
$f:E\to F$.
\end{prop}

Proof. Pick $f:E\to F$ such that $J\cong \im f$.
Composing a surjection $E'\to E$ and an injection
$F\into F'$, we get another map
 $f':E'\to F'$ such that $J\cong \im f'$.
We get two descent problems, ${\mathbf D}$ and ${\mathbf D}'$.
We claim that
$$
C^*\bigl({\rm Sci}({\mathbf D}), F\bigr)=
C^*\bigl({\rm Sci}({\mathbf D}'), F'\bigr).
$$
This follows from (\ref{main.thm.v1}) and the obvious maps
$$
C^*(X, E')\onto C^*(X, E)\to C^*(X, F)\into C^*(X, F'). \qed
$$

\subsection*{Proof of Theorem \ref{main.thm.v0}}{\ }

In order to get an idea of the proof, assume first that
$X,Y$ are normal and let $Y\to W\to X$ denote the Stein factorization.
We  first study which sections over $Y$ descend to
$W$ and then try to descend them to $X$.

If we look over a single point $w\in W$, the the question 
is answered by (\ref{wronsk.lem}).
Working in our family, this means passing from
$Y\to W$ to the $(n+1)$-fold fiber product
$Y\times_W\times\cdots\times_W Y$.
The fiber product can be rather singular in general,
so this will work only over a dense  open subset of $W$.

Going from $W$ to $X$ is easy if we work locally analytically.
In this case $W\to X$ is a local isomorphism over an
 open subset of $W$, thus every question over $W$ can be rewritten
as a question over $X$. This will not work well algebraically, but
there are no problems if $W\to X$ is Galois.

The point of (\ref{red.step.main.2}) is to show that
by passing to a suitable scion, the above considerations apply,
at least over a dense open subset of $X$.

Then we finish by a straightforward dimension induction
(\ref{pf.of.main.thm}).

\begin{prop}\label{red.step.main.2}
 Let ${\mathbf D}=\bigl(p:Y\to X, f:p^*E\to F\bigr)$
 be a  descent problem.
Then there is a  closed
 algebraic subvariety $Z\subset X$ with $\dim Z<\dim X$ and 
 a  scion 
$\tilde {\mathbf D}=\bigl(\tilde p:\tilde Y\to X, \tilde f:
\tilde p^*E\to \tilde F\bigr)$
with surjective structure map $\tilde r:\tilde Y\to Y$
and with the following properties.

Let $X=\cup_{i\in I} X_i$ be the irreducible components.
For every $i\in I$ let  $\tilde Y_i\subset \tilde Y$ be  the closure of
$\tilde p^{-1}(X_i\setminus Z)$ and
$\tilde {\mathbf D}_i$  the restriction of $\tilde {\mathbf D}$ to 
$\tilde Y_i$.
Then, for  every $i\in I$,
\begin{enumerate}
\item a  finite group $G_i$ acts on  $\tilde {\mathbf D}_i$,
\item there is a $G_i$-equivariant factorization 
 $\tilde p_i:\tilde Y_i\stackrel{\tilde q_i}{\to} \tilde W_i
\stackrel{\tilde w_i}{\to} X_i$,
\item  over $X_i\setminus Z$, the map $\tilde w_i: \tilde W_i\to X_i$ is  finite
and Galois with group $G_i$,
\item there is a $G_i$-equivariant quotient bundle
$\tilde w_i^*E\to \tilde E_i$ such that $\tilde f_i$ factors as
$\tilde p_i^*E\onto \tilde q_i^*\tilde E_i\cong \tilde F_i$.
\end{enumerate}
\end{prop}

Proof. 
We may harmlessly assume that $p(Y)$ is  dense in $X$.

After we construct  $\tilde {\mathbf D}$,
the plan 
is to make sure that $Z$ contains all of its  ``singular'' points. 
In the original setting
$Z$ is the set where the map
$(f_1,\dots, f_r):\o_X^r\to \o_X$ has rank 0.
In the general case, we need to include points over which
$\tilde f$ drops rank and also points over which $\tilde p$ drops rank.
During the proof we gradually add more and more
irreducible components to $Z$ as needed.

{\it Step 0.} To start with, we  add to $Z$  the locus where $X$ is not normal
and    the $p(Y_j)$ where
$Y_j\subset Y$ is an irreducible component
that does not dominate any of the 
irreducible components  of $X$.
In the conclusions, the different $\tilde{\mathbf D}_i$ have
no effect on each other, hence we can work with them
one at a time. 
We construct each $\tilde{\mathbf D}_i$ separately, and then let
$\tilde{\mathbf D}$ be the disjoint union of the
$\tilde{\mathbf D}_i$ for $i\in I$ and of ${\mathbf D}|_Z$. 

For simplicity of notation, we drop the index $i$.
We  thus assume that $X$ is irreducible and
every  irreducible component of $Y$ dominates $X$.
We may assume that $Y$ is normal,
take  the Stein factorization $p:Y\stackrel{q}{\to} W\stackrel{s}{\to} X$
and set $m=\deg(W/X)$.
In several steps we construct
the following diagram
$$
\begin{array}{rclclclcc}
\bigl(\tilde q^{(m)}\bigr)^*\bar E^{(m)} & \cong & \tilde F^{(m)} &&
\bar F^{(m)} &  & F^{(m)} &  & F\\
&&  \ \downarrow  &&\ \downarrow  &&\ \downarrow  &&\downarrow  \\
\bigl(t^{(m)}\circ s^{(m)}\bigr)^*E\onto \bar E^{(m)}  & & \tilde Y^{(m)}_X
& & \bar Y^{(m)}_X & 
\stackrel{t_Y^{(m)}}{\to} & Y^{(m)}_X & \stackrel{\pi_i^{(m)}}{\to} & Y\\
&\searrow & \ \downarrow\tilde q^{(m)}
  &&\ \downarrow \bar q^{(m)}  &&\ \downarrow  q^{(m)}&&\hphantom{p}\downarrow p  \\
&&\tilde W^{(m)}_X & =  & \bar W^{(m)}_X & 
\stackrel{t^{(m)}}{\to} & W^{(m)}_X &  
\stackrel{s^{(m)}}{\to}  & X
\end{array}
$$

{\it Step 1: Constructing $W^{(m)}_X$ and its column.}

Let $s:W\to X$ be a finite morphism of (possibly reducible) varieties.

Consider the $m$-fold fiber product
$W^m_X:=W\times_X\cdots\times_XW$ with coordinate projections
$\pi_i:W^m_X\to W$. For every $i\neq j$,
let $\Delta_{ij}\subset W^m_X$ be the preimage of the diagonal
 $\Delta\subset W\times_XW$ under the map $(\pi_i,\pi_j)$.
Let $W^{(m)}_X\subset W^m_X$ be the 
union of the dominant components in the 
closure of
$W^m_X\setminus\cup_{i\neq j}\Delta_{ij}$
with projection $s^{(m)}:W^{(m)}_X\to X$.
The symmetric group $S_m$ acts on $W^{(m)}_X$ by permuting the factors.

 Let
 $X^0\subset X$  be the largest Zariski open subset over which $s$ is smooth.
If $x\in X^0$ then
$\bigl(s^{(m)}\bigr)^{-1}(x)$ consists of
ordered $m$-element subsets of $s^{-1}(x)$,
thus 
$S_m$ acts transitively on $\bigl(s^{(m)}\bigr)^{-1}(x)$ if
$|s^{-1}(x)|=m$.

Let now $p:Y\to X$ be as above
with Stein factorization $p:Y\stackrel{q}{\to} W\stackrel{s}{\to} X$.
Let $Y^m_X$ denote the $m$-fold fiber product $Y\times_X\cdots\times_XY$
 with coordinate projections
$\pi_i:Y^m_X\to Y$.

Let  $Y^{(m)}_X\subset Y^m_X$ 
denote the  dominant parts of the 
preimage of  $W^{(m)}_X$ under the natural map
$q^m:Y^m_X\to W^m_X$ with projection $p^{(m)}:Y^{(m)}_X\to X$.
Note that, for general $x\in X$, 
$S_m$ acts transitively on the  irreducible components of 
$\bigl(p^{(m)}\bigr)^{-1}(x)$.

Let $F$ be a locally free sheaf on $Y$.
Then $\oplus_i \pi_i^*F$ is a   locally free sheaf on  $Y^m_X$.
Its restriction to  $Y^{(m)}_X$ is denoted by $F^{(m)}$. 

The $S_m$-action on  $Y^{(m)}_X$ naturally lifts to
an  $S_m$-action on  $F^{(m)}$. 
From $f:p^*E\to F$  we get an $S_m$-invariant
 map  of locally free sheaves
$f^{(m)}:\bigl(p^{(m)}\bigr)^*E\to F^{(m)}$. 
For each $m$ we get a scion of  ${\mathbf D}$
$$
{\mathbf D}^{(m)}:=\bigl(p^{(m)}:Y^{(m)}_X\to X, 
f^{(m)}:\bigl(p^{(m)}\bigr)^*E\to F^{(m)}\bigr).
$$
\medskip

{\it Step 2: Constructing $\bar W^{(m)}_X$ and its column.}

More generally, let 
${\mathbf D}=\bigl(q:Y\to W, f:q^*E\to F\bigr)$ be  a descent problem.
(Note that the base is $W$ instead of $X$.)
Assume that $W$ is irreducible. 
Consider the coherent sheaf
$E':= \im\bigl[E\to  q_* F\bigr]$.

Let $Gr(d,E)\to W$ be the universal Grassmann bundle of 
rank $d$ quotients of $E$ where $d$ is the  rank of $E'$ at a general point. 
At a general point $w\in W$, $E(w)\onto  E'(w)$
 is such a quotient.
Thus $E'$ gives a rational map
$W\map Gr(d,E)$, defined on a  dense 
 open subset. Let $\bar W\subset Gr(d,E)$
 denote the
closure of its image  and $t:\bar W\to W$
 the projection.
Then  $t$ is a proper birational morphism
and we have a decomposition
 $$
t^*q: t^*E\stackrel{s}{\onto} \bar E\stackrel{j}{\into}  t^*E'
$$
where $\bar E$  is a locally free sheaf of rank $d$ on $\bar W$,
$s$ is a rank $d$ surjection everywhere and $j$
is a rank $d$ injection on a  dense  open subscheme.

Applying this to ${\mathbf D}^{(m)}$, with $W^{(m)}_X$
playing the role  of the base, 
 we obtain
$\bar {\mathbf D}^{(m)}$.
\medskip

{\it Step 3: Constructing $\tilde W^{(m)}_X$ and its column.}

More generally, let
${\mathbf D}=\bigl(q:Y\to W, f:q^*E\to F\bigr)$ be  a descent problem.
Assume that $W$ and  the generic fiber of $q$
are irreducible and  $E\to q_*F$
is an injection. We construct a  scion
$$
\tilde {\mathbf D}=\bigl(\tilde q:\tilde Y\to W, \tilde f:
\tilde q^*E\to \tilde F\bigr)
$$
with surjective structure  map
such that $\tilde f$  is an isomorphism.

 Set $n=\rank E$ and 
let $Y^{n+1}_W$ be the  union of the dominant components 
 of the
$n+1$-fold fiber product of $Y\to W$ with
coordinate projections $\pi_i$. 
Let $\tilde q:Y^{n+1}_W\to W$ be the map given by any of the $q\circ \pi_i$.
Consider the diagonal map
$$
\tilde f: \tilde q^*E\to \tsum_{i=1}^{n+1} \pi_i^*F
$$
which is an injection 
over a  dense  open set $Y^0\subset Y^{n+1}_W$ by assumption.
Using  (\ref{rel.desc.prob.sefn}.3)
we can  replace
$\tsum_{i=1}^{n+1} \pi_i^*F$ by $\tilde q^*E$.

Applying this to $\bar {\mathbf D}^{(m)}$ we obtain
$\tilde {\mathbf D}^{(m)}$.
\qed

\begin{prop}\label{red.step.main.1}
 Let 
$\tilde {\mathbf D}=\bigl(\tilde p:\tilde Y\to X, \tilde f:
\tilde p^*E\to \tilde F\bigr)$
 be a  descent problem and $Z\subset X$ a 
 closed algebraic subvariety.
Let $X=\cup_i X_i$ be the irreducible components
and assume that $X_i\cap X_j\subset Z$ for every $i\neq j$.
Let $\tilde Y_i\subset Y$ be the closure of
$\tilde p^{-1}(X_i\setminus Z)$ and
$\tilde {\mathbf D}_i$  the restriction of $\tilde {\mathbf D}$
 to $\tilde Y_i$.
Assume that for every $i$ 
(\ref{red.step.main.2}.1--4) hold.

Then $\tilde{\mathbf D}$ is finitely determined relative to $Z$.
\end{prop}

Proof. Let $\Psi_Y\in C^*(Y, F)$ be a section  that
 vanishes on  $p^{-1}( Z)$ such that 
(\ref{findet.defn}.2) holds.
We can uniquely write $\Psi_Y=\sum \Psi_i$
where $\supp \Psi_i\subset \tilde Y_i$.
It is thus enough to write $\Psi_i=f\circ p^*(\psi_{i,X})$
for each $i$. For a fixed $i$, we need to do this
over $X_i$ and then extend $\psi_{i,X}$ to $X$ by setting it
zero on the complement. Thus it is sufficient to work with
 one $\tilde {\mathbf D}_i$ at a time.

Using the isomorphism $\tilde q_i^*\tilde E_i\cong \tilde F_i$,
  $\Psi_i$ can be identified with a section 
 $\tilde\Psi_i$ of  $\tilde q_i^*\tilde E_i$.
The conditions (\ref{findet.defn}.2) now imply that
 $\tilde\Psi_i$
 is constant on the fibers of $\tilde Y_i\to \tilde W_i$
and is  $G_i$-invariant.
Thus $\tilde\Psi_i$ is the pull-back of a  $G_i$-invariant 
section $\tilde\Psi_{W,i}$ of
$\tilde E_i$ that vanishes on the preimage of $Z$. 
Using a  $G_i$-invariant $C^*$-splitting of 
$\tilde w_i^*E\onto \tilde E_i$
we can think of $\tilde\Psi_{W,i}$ as a  $G_i$-invariant section of
$\tilde w_i^*E$. Therefore   $\tilde\Psi_{W,i}$ descends to a   section
$\psi_{X,i} \in C^*\bigl(X_i, E\bigr)$
that vanishes on $Z$. \qed

\begin{say}[Proof of Theorem \ref{main.thm.v0}]\label{pf.of.main.thm}

We use induction on the dimension of $X$. 
If $\dim X=0$ then we are done by (\ref{wronsk.lem}).

In  general, construct $Z$ and 
$\tilde{\mathbf D}$ as in (\ref{red.step.main.2}). 
Let  $\tilde{\mathbf D}_Z$ denote the restriction of  
$\tilde{\mathbf D}$ to $Z$.
By induction, it has a finitely determined scion
whose structure map  is surjective;  we denote it
$\bigl(\tilde{\mathbf D}_Z\bigr)^{\sim}$. Let ${\mathbf D}_s$
be the disjoint union of  $\tilde{\mathbf D}$ and of
$\bigl(\tilde{\mathbf D}_Z\bigr)^{\sim}$.

Pick $\Phi_s\in C^*\bigl(Y_s, F_s)$ and assume that
it satisfies the conditions (\ref{findet.defn}.2).
Its restriction to $\bigl(\tilde{\mathbf D}_Z\bigr)^{\sim}$
also  satisfies the conditions (\ref{findet.defn}.2), hence
there is a section $\phi_Z\in C^*\bigl(Z, E|_{Z}\bigr)$
whose pull-back to $p_s^{-1}(Z)$ equals the restriction
of $\Phi_s$. 
(A priori this holds only over $\bigl(\tilde Y_Z\bigr)^{\sim}$,
but since the structure map
$\bigl(\tilde Y_Z\bigr)^{\sim}\to \tilde Y_Z$ is surjective,
it also holds over $\tilde Y_Z$.)

By (\ref{basic.props}), we can lift $\phi_Z$ to a section
$\phi_X \in C^*\bigl(X, E\bigr)$.
Consider next  
$$
\Psi_s:=\Phi_s-f_s\bigl(p_s^*\phi_X\bigr)\in  
C^*\bigl(Y_s, F_s\bigr).
$$  
By construction, it vanishes along $p_s^{-1}(Z)$.
By (\ref{red.step.main.1}),
  $\tilde{\mathbf D}$ is finitely determined relative to $Z$,
hence  we can write $\Psi_s=f_s\circ p_s^*(\psi_X)$
for some $\psi_X \in C^*\bigl(X, E\bigr)$.
Thus $\Phi_s=f_s\circ p_s^*\bigl(\phi_X+\psi_X\bigr)$.
 \qed
\end{say}

\begin{ack} I learned about continuous closures from a 
very interesting lecture of
 M.~Hochster. 
The questions and  comments of H.~Brenner,
 Ch.~Fefferman, A.~N\'emethi and T.~Szamuely
were very helpful in formulating many of the definitions and results.
Partial financial support  was provided by  the NSF under grant number 
DMS-0758275.
\end{ack}

\bibliography{crefs}

\providecommand{\bysame}{\leavevmode\hbox to3em{\hrulefill}\thinspace}
\providecommand{\MR}{\relax\ifhmode\unskip\space\fi MR }
\providecommand{\MRhref}[2]{%
  \href{http://www.ams.org/mathscinet-getitem?mr=#1}{#2}
}
\providecommand{\href}[2]{#2}
\begin{thebibliography}{BCR98}

\bibitem[BCR98]{bcr}
Jacek Bochnak, Michel Coste, and Marie-Fran{\cedilla{c}}oise Roy, \emph{Real
  algebraic geometry}, Ergebnisse der Mathematik und ihrer Grenzgebiete,
  vol.~36, Springer-Verlag, Berlin, 1998, Translated from the 1987 French
  original, Revised by the authors. \MR{MR1659509 (2000a:14067)}

\bibitem[Bre06]{brenner}
Holger Brenner, \emph{Continuous solutions to algebraic forcing equations},
  http://www.citebase.org/abstract?id=oai:arXiv.org:0608611, 2006.

\bibitem[FK10]{fef-kol}
Charles Fefferman and J{\'a}nos Koll{\'a}r, \emph{Continuous linear
  combinations of polynomials}, (in preparation), 2010.

\bibitem[Har77]{hartsh}
Robin Hartshorne, \emph{Algebraic geometry}, Springer-Verlag, New York, 1977,
  Graduate Texts in Mathematics, No. 52. \MR{0463157 (57 \#3116)}

\bibitem[Hoc10]{hochster}
Melvin Hochster, \emph{(personal communication)}, 2010.

\bibitem[Kol96]{rc-book}
J{\'a}nos Koll{\'a}r, \emph{Rational curves on algebraic varieties}, Ergebnisse
  der Mathematik und ihrer Grenzgebiete. 3. Folge., vol.~32, Springer-Verlag,
  Berlin, 1996. \MR{1440180 (98c:14001)}

\bibitem[Kol10]{k-obw}
J{\'a}nos Koll{\'a}r, \emph{Continuous closure of sheaves}, Mathematisches
  Forschungsinstitut Oberwolfach, Report No. 27, 2010, pp.~41--43.

\end{thebibliography}

\vskip1cm

\noindent Princeton University, Princeton NJ 08544-1000

{\begin{verbatim}kollar@math.princeton.edu\end{verbatim}}

\end{document}